\documentclass[12pt,leqno,draft]{article}

\def\radius{\mathop{\rm radius}}

\newtheorem{theorem}{Theorem}
\newtheorem{lemma}[theorem]{Lemma}
\newtheorem{proposition}[theorem]{Proposition}
\newtheorem{sublemma}[theorem]{Sublemma}
\newtheorem{definition}[theorem]{Definition}
\newtheorem{corollary}[theorem]{Corollary}
\newtheorem{problem}[theorem]{Problem}
\newtheorem{remark}[theorem]{Remark}
\newtheorem{claim}[theorem]{Claim}
\newtheorem{assumptions}[theorem]{Assumptions}
\newtheorem{examples}[theorem]{Examples}
\newtheorem{question}[theorem]{Question}
\newtheorem{sassumptions}[theorem]{Standing Assumptions}
\newtheorem{sassumption}[theorem]{Standing Assumption}
\newtheorem{conjecture}[theorem]{Conjecture}

\newcommand{\begintheorem}{\addtocounter{equation}{1}\begin{theorem}}
\newcommand{\beginlemma}{\addtocounter{equation}{1}\begin{lemma}}
\newcommand{\beginproposition}{\addtocounter{equation}{1}\begin{proposition}}
\newcommand{\beginsublemma}{\addtocounter{equation}{1}\begin{sublemma}}
\newcommand{\begindefinition}{\addtocounter{equation}{1}\begin{definition}}
\newcommand{\begincorollary}{\addtocounter{equation}{1}\begin{corollary}}
\newcommand{\beginproblem}{\addtocounter{equation}{1}\begin{problem}}
\newcommand{\beginremark}{\addtocounter{equation}{1}\begin{remark}}
\newcommand{\beginclaim}{\addtocounter{equation}{1}\begin{claim}}
\newcommand{\beginassumptions}{\addtocounter{equation}{1}\begin{assumptions}}
\newcommand{\beginexamples}{\addtocounter{equation}{1}\begin{examples}}
\newcommand{\beginquestion}{\addtocounter{equation}{1}\begin{question}}
\newcommand{\beginsassumptions}{\addtocounter{equation}{1}\begin{sassumptions}}
\newcommand{\beginsassumption}{\addtocounter{equation}{1}\begin{sassumption}}
\newcommand{\beginconjecture}{\addtocounter{equation}{1}\begin{conjecture}}

\begin{document}

\title{Some questions related to fractals}

\author{Stephen Semmes}

\date{}

\maketitle

\subsubsection*{Classes of fairly smooth functions}

	A fundamental aspect of harmonic analysis on Euclidean spaces
deals with various classes of real or complex-valued functions, such
as the $C^k$ classes of functions which are continuous and
continuously differentiable up to order $k$, where $k$ is a
nonnegative integer.  If $k = 0$, this simply means that the function
is continuous.  Of course continuity of a function makes sense on any
metric space, or on topological spaces more generally, while the
notion of \emph{derivatives} entails more structure.  Compare with
\cite{Rudin1}.

	We can be more precise and consider $C^{k, \alpha}$ classes of
functions, where $k$ is a nonnegative integer and $\alpha$ is a real
number, $0 \le \alpha \le 1$.  If $\alpha = 0$, then $C^{k, \alpha}$
is taken to mean the same as $C^k$; when $\alpha > 0$, $C^{k, \alpha}$
consists of the $C^k$ functions with the extra property that their
$k$th order derivatives are locally \emph{H\"older continuous of order
$\alpha$}.  For this let us recall that a function $h(x)$ on a metric
space $(M, d(x,y))$ is H\"older continuous of order $\alpha$ if there
is a nonnegative real number $A$ such that
\begin{equation}
\label{Holder continuity of order alpha}
	|h(x) - h(y)| \le A \, d(x,y)^\alpha \quad \hbox{for all } x, y \in M 
\end{equation}
This makes sense on any metric space, so that $C^{0, \alpha}$ classes
of functions make sense on any metric space, but $C^{k, \alpha}$
involves more structure when $k \ge 1$.

	The condition (\ref{Holder continuity of order alpha}) also
makes sense for $\alpha > 1$, but on a Euclidean space it would then
imply that the function is constant, because the first derivatives
would be identically equal to $0$.  One can also check this fact more
directly.  On a general metric space, this is not true in general.
However, it is true under fairly mild conditions, e.g., if the metric
space is connected and there are enough rectifiable curves around,
because such a function would have to be constant on each one of them.
Cantor sets and snowflake curves are basic examples of metric spaces
in which there are a lot of functions which satisfy (\ref{Holder
continuity of order alpha}) with $\alpha > 1$.

	Another smoothness class for functions on ${\bf R}^n$
is given by the \emph{Zygmund condition}
\begin{equation}
\label{Zygmund condition}
	|f(x + v) + f(x - v) - 2 f(x)| \le L \, |v| 
				\quad \hbox{for all } x, v \in {\bf R}^n,
\end{equation}
where $f$ is a continuous function on ${\bf R}^n$ and $L$ is a
nonnegative real number.  This condition is implied by (\ref{Holder
continuity of order alpha}) when $\alpha = 1$, with $L = 2 A$, but the
converse does not work, even locally.  It can be shown that if a
continuous function $f$ satisfies (\ref{Zygmund condition}), then it
is locally H\"older continuous of order $\alpha$ for each $\alpha <
1$.

	Note that the Zygmund condition is not defined on general
metric spaces.  Here is a nice reformulation of it, although not
exactly with the same constant, in terms of affine functions on ${\bf
R}^n$.  A function $f$ on ${\bf R}^n$ satisfies the Zygmund condition
for some $L$ if and only if there is a nonnegative real number $L'$
such that for every ball $B$ in ${\bf R}^n$ there is an affine
function $a(x)$ such that
\begin{equation}
\label{Zygmund condition, revisited}
	\sup_{x \in B} |f(x) - a(x)| \le L' \radius(B).
\end{equation}
It is very easy to go from this to (\ref{Zygmund condition}) with $L = L'$,
but the other direction is more complicated.

	If we replace $|v|$, $\radius(B)$ on the right sides of
(\ref{Zygmund condition}), (\ref{Zygmund condition, revisited}) with
$|v|^{1 + \alpha}$, $\radius(B)^{1 + \alpha}$, respectively, where $0
< \alpha \le 1$, then the corresponding conditions can be shown to be
equivalent to $f$ being continuously differentiable with first
derivatives H\"older continuous of order $\alpha$.  To go from
H\"older continuous first derivatives to the other conditions is
basically a matter of calculus, but the other direction is again more
complicated.  If one replaces $|v|$, $\radius(B)$ on the right sides
of (\ref{Zygmund condition}), (\ref{Zygmund condition, revisited})
with $|v|^{1 + \alpha}$, $\radius(B)^{1 + \alpha}$ when $\alpha > 1$,
then the resulting conditions imply that the second derivatives of $f$
are equal to $0$, and $f$ is affine.  To accommodate higher powers of
$|v|$ or $\radius(B)$, one can use differences of higher order in
place of the second difference in (\ref{Zygmund condition}), or
polynomials of higher degree in place of affine functions in
(\ref{Zygmund condition, revisited}).

	Some basic references concerning these matters are \cite{Duren,
Krantz2, MS, SS1, SS3, St1, SW, TW, Zygmund}.

	All of this uses the special structure of Euclidean spaces, to
define differences of order at least $2$, or to have affine functions
and polynomials of higher degree, rather than just constants.  What
about other metric spaces, and the structure that they might have?
This leads to a lot of questions.

	In the case of Heisenberg groups, other nilpotent Lie groups,
and sub-Riemannian spaces more generally, a lot of theory has been
developed along these lines.  There are various smoothness classes
analogous to classical ones on Euclidean spaces which are adapted to
the new geometry.  A basic notion is that $C^\infty$ functions or
polynomials in general might be the same as before, but more precise
degrees of smoothness are now different, and different kinds of
degrees of polynomials are used.  See \cite{FS, Krantz3}, for instance.

	What about fractal sets sitting inside of Euclidean spaces?
For these the usual affine functions and polynomials can be used,
just restricted to the fractal set under consideration.  One
can also define classes like $C^{k, \alpha}$ on the set by taking
restrictions of $C^{k, \alpha}$ functions from the ambient Euclidean
space.  Well-known extension theorems of Whitney characterize
functions on a closed set in a Euclidean space which can be extended
to $C^{k, \alpha}$ functions on the whole Euclidean space, as in
\cite{K-P1, St1}.

	Smoothness inherited from the ambient space in this manner
is quite different from the kind of regularity discussed in
\cite{Kigami, Strichartz2}, for instance.

	Other very interesting examples with their own special behavior
can be found in \cite{Bourdon, BP, Laakso1, Laakso2}.

\subsubsection*{Semi-Markovian spaces}

	Gromov \cite{misha-hypgroups, C-P} introduced a remarkable
notion of ``semi-Markovian spaces'' for which there are combinatorial
``presentations'', and which includes a number of standard fractals.
This appears to have a lot of room for interesting analysis of
functions.

\subsubsection*{Families of fractals}

	It seems to me that there are a lot of issues around the
general theme of ``families of fractals''.  A basic set-up would be to
have two metric spaces $T$ and $P$ and a mapping $\pi : T \to P$,
where the fibers $\pi^{-1}(p)$, $p \in P$, would be the fractals in
the family.  It is customary to call $T$ the ``total space'' for the
family of fractals, and $P$ the ``space of parameters''.  One might
ask that $\pi$ be Lipschitz, which is to say H\"older continuous of
oder $1$.  One could also ask for some nondegeneracy conditions, along
the lines of the nonlinear quotient mappings discussed in \cite{BJLPS,
BL, JLPS1, JLPS2}, or $(\tau, \rho)$-regular mappings or noncollapsing
mappings as in \cite{D-S2}.  For $(\tau, \rho)$-regular mappings,
$\tau$, $\rho$ could be some kind of dimensions for $T$, $P$,
respectively.

	Here is a somewhat more specific type of situation.  Suppose
that $n$ is a positive integer, and consider the Cartesian product
${\bf R}^n \times {\bf R}$.  Let $\lambda : {\bf R}^n \times {\bf R}
\to {\bf R}$ be the standard projection onto the last coordinate, so
that $\lambda(x, t) = t$ for all $x$ in ${\bf R}^n$ and $t$ in ${\bf
R}$.  Let $E$ be a subset of ${\bf R}^n \times {\bf R}$, and assume
for instance that $E$ is compact and that $\lambda(E) = [0,1]$.  The
sets $E(t) = \lambda^{-1}(t)$, $0 \le t \le 1$, would be the fractals
in the family.  Let $d$ be a positive real number, and suppose that
$E(t)$ has Hausdorff dimension $d$ for all $t$ in $[0,1]$, or almost
all $t$.  Then $E$ itself should have Hausdorff dimension at least $d
+ 1$, and if the Hausdorff dimension of $E$ is equal to $d + 1$,
then that reflects a kind of regularity in the situation.  Of course
there are a number of variants of this.

\end{document}